\documentclass[leqno,12pt]{article}
\usepackage{amsmath, amssymb,soul}
\usepackage{amsthm} 

\theoremstyle{plain} 
\newtheorem{theorem}{\indent\sc Theorem}[section]
\newtheorem{lemma}[theorem]{\indent\sc Lemma}

\newtheorem{proposition}[theorem]{\indent\sc Proposition}

\numberwithin{equation}{section}
\theoremstyle{definition} 
\newtheorem{definition}[theorem]{\indent\sc Definition}
\newtheorem{remark}[theorem]{\indent\sc Remark}
\newtheorem{example}[theorem]{\indent\sc Example}


\def\({\left( }
\def\){\right)}
\def\<{\left< }
\def\>{\right>}

\title{Warped product pointwise bi-slant submanifolds in metallic Riemannian manifolds}
\author{Cristina E. Hretcanu and Adara M. Blaga}
\date{}

\begin{document}

\maketitle

\markboth{{\small\it {\hspace{2cm} Warped product pointwise bi-slant submanifolds in metallic Riemannian manifolds}}}{\small\it{Warped product pointwise bi-slant submanifolds in metallic Riemannian manifolds
\hspace{2cm}}}

\textbf{Abstract:} In this paper, we study some properties of warped product pointwise bi-slant submanifolds in locally metallic Riemannian manifolds and we construct some examples in Euclidean spaces.\\

{ 
\textbf{2020 Mathematics Subject Classification:}
53B20, 53B25, 53C42, 53C15.
}

{ 
\textbf{Keywords:}
Metallic Riemannian structure; Warped product bi-slant submanifold; Pointwise slant submanifold.
}


\section{Introduction}

Metallic Riemannian manifolds and their submanifolds were defined and investigated by C. E. Hretcanu, M. Crasmareanu and A. M. Blaga in (\cite{Hr5}, \cite{Hr4}), as a generalization of Golden Riemannian manifolds studied in (\cite{CrHr}, \cite{Hr3}, \cite{Hr2}). The authors of the present paper obtained some properties of invariant, anti-invariant and slant submanifolds (\cite{Blaga_Hr}), semi-slant submanifolds (\cite{Hr6}) and, respectively, hemi-slant submanifolds (\cite{Hr7}) in metallic and Golden Riemannian manifolds and they provided some integrability conditions for the distributions involved in these types of submanifolds. Moreover, properties of metallic and Golden warped product Riemannian manifolds were presented in some previous works of the authors (\cite{Blaga1}, \cite{Blaga2}, \cite{Hr10}). In the last years, the study of submanifolds in metallic Riemannian manifolds has been continued by many authors (\cite{bea1}, \cite{fee1}, \cite{fee2}), which introduced the notion of lightlike submanifold of a metallic semi-Riemannian manifold.

\section{Preliminaries}

The name of metallic number is given to the positive solution of the equation
\linebreak
$x^{2}-px-q=0$, which is $\sigma _{p,q}=\frac{p+\sqrt{p^{2}+4q}}{2}$ (\cite{Spinadel}), where $p$ and $q$ are positive integer values. The metallic structure is a particular case of polynomial structure on a manifold, which was generally defined in (\cite{Goldberg2}, \cite{Goldberg1}).

Let $\overline{M}$ be an $m$-dimensional manifold endowed with a tensor field $J$ of type $(1,1)$. Then $J$ is called a \textit{metallic structure} if it satisfies:
\begin{equation}\label{e1}
J^{2}= pJ+qI,
\end{equation}
for $p$, $q\in\mathbb{N}^*$, where $I$ is the identity operator on $\Gamma(T\overline{M})$. If a Riemannian metric $\overline{g}$ is $J$-compatible, i.e.:
\begin{equation} \label{e2}
\overline{g}(JX, Y)= \overline{g}(X, JY),
\end{equation}
for any $X$, $Y \in \Gamma(T\overline{M})$, then $(\overline{M},\overline{g},J)$ is called a {\it metallic Riemannian manifold} (\cite{Hr4}).
In this case, $\overline{g}$ verifies:
\begin{equation} \label{e3}
\overline{g}(JX, JY)=\overline{g}(J^{2}X, Y) =p \overline{g}(JX,Y)+q \overline{g}(X,Y),
\end{equation}
for any $X$, $Y \in \Gamma(T\overline{M})$.
\normalfont

For $p=q=1$ one obtain the \textit{Golden structure} $J$ which satisfies $J^{2}= J + I$. If ($\overline{M}, \overline{g})$ is a Riemannian manifold endowed with a Golden structure $J$ such that the Riemannian metric $\overline{g}$ is $J$-compatible, then $(\overline{M},\overline{g},J)$ is called a {\it Golden Riemannian manifold} (\cite{CrHr}).

Let $M$ be an isometrically immersed submanifold in the metallic Riemannian manifold ($\overline{M}, \overline{g},J)$.
The tangent space $T_x\overline{M}$ of $\overline{M}$ in a point $x \in M$ can be decomposed into the direct sum
$T_x\overline{M}=T_x M\oplus T_x^{\perp}M,$ for any $x\in M$, where $T_{x}^{\bot }M$ is the normal space of $M$ in $x$.
Let $i_{*}$ be the differential of the immersion $i: M \rightarrow\overline{M}$. Then the induced Riemannian metric $g$ on $M$ is given by $g(X, Y)=\overline{g}(i_{*}X, i_{*}Y)$, for any $X$, $Y \in \Gamma(TM)$. In all the rest of the paper, we shall denote by $X$ the vector field $i_{*}X$, for any $X \in \Gamma(TM)$.

For any $X \in \Gamma(TM)$, let $TX:=(J X)^T$ and $NX:=(J X)^{\perp}$ be the tangential and normal components, respectively, of $JX$ and for any $V \in \Gamma(T^{\perp}M)$, let $tV:=(J V)^{T}$, $nV:=(J V)^{\perp}$ be the tangential and normal components, respectively, of $JV$. Then we have:
\begin{equation}\label{e4}
JX = TX + NX,
\end{equation}
\begin{equation}\label{e5}
JV = tV + nV,
\end{equation}
for any $X \in \Gamma(TM)$ and $V \in \Gamma(T^{\perp}M)$.

The maps $T$ and $n$ are $\overline{g}$-symmetric (\cite{Blaga_Hr}):
\begin{equation}\label{e6}
\overline{g}(TX,Y)=\overline{g}(X,TY),
\end{equation}
\begin{equation}\label{e7}
\overline{g}(nU,V)=\overline{g}(U,nV)
\end{equation}
and
\begin{equation}\label{e8}
 \overline{g}(NX,V)=\overline{g}(X,tV),
\end{equation}
for any $X$, $Y \in \Gamma(TM)$ and $U$, $V\in \Gamma(T^{\perp}M)$.

\pagebreak
We also obtain (\cite{Hr7}):
\begin{equation} \label{e9}
T^{2}X = pTX+qX-tNX,
\end{equation}
\begin{equation} \label{e10}
pNX= NTX+nNX,
\end{equation}
\begin{equation} \label{e11}
n^{2}V =  pnV+qV-NtV,
\end{equation}
\begin{equation} \label{e12}
ptV= TtV+tnV,
\end{equation}
for any $X \in \Gamma(TM)$ and $V \in \Gamma(T^{\bot}M)$.

Let $\overline{\nabla}$ and $\nabla $ be the Levi-Civita connections on $(\overline{M},\overline{g})$ and on its submanifold $(M,g)$, respectively. The Gauss and Weingarten formulas are given by:
\begin{equation}\label{e13}
\overline{\nabla}_{X}Y=\nabla_{X}Y+h(X,Y),
\end{equation}
\begin{equation}\label{e14}
\overline{\nabla}_{X}V=-A_{V}X+\nabla_{X}^{\bot}V,
\end{equation}
for any $X$, $Y \in \Gamma(TM)$ and $V \in \Gamma(T^{\bot}M)$, where $h$ is the second fundamental form and $A_{V}$ is the shape operator, which satisfy
\begin{equation}\label{e15}
 \overline{g}(h(X, Y),V)=\overline{g}(A_{V}X, Y).
 \end{equation}

For any $X$, $Y \in \Gamma(TM)$, the covariant derivatives of $T$ and $N$ are given by:
\begin{equation}\label{e16}
(\nabla_{X}T)Y=\nabla_{X}TY - T(\nabla_{X}Y),
\end{equation}
\begin{equation}\label{e17}
(\overline{\nabla}_{X}N)Y=\nabla_{X}^{\bot}NY - N(\nabla_{X}Y).
\end{equation}

For any $X \in \Gamma(TM)$ and $V \in \Gamma(T^{\bot}M)$, the covariant derivatives of $t$ and $n$ are given by:
\begin{equation}\label{e18}
(\nabla_{X}t)V=\nabla_{X}tV - t(\nabla_{X}^{\bot}V),
\end{equation}
\begin{equation}\label{e19}
(\overline{\nabla}_{X}n)V=\nabla_{X}^{\bot}nV - n(\nabla_{X}^{\bot}V).
\end{equation}

From $(\ref{e1})$ we obtain:
\begin{equation} \label{e20}
\overline{g}((\overline{\nabla}_XJ)Y,Z)=\overline{g}(Y,(\overline{\nabla}_XJ)Z),
\end{equation}
for any $X$, $Y$, $Z\in \Gamma(T\overline{M})$, which implies (\cite{Blaga3}):

\pagebreak
\begin{equation}\label{e21}
\overline{g}((\nabla_X T)Y,Z)=\overline{g}(Y,(\nabla_X T)Z),
\end{equation}
\begin{equation}\label{e22}
\overline{g}((\overline{\nabla}_{X}N)Y,V )= \overline{g}(Y,(\nabla_{X}t)V),
\end{equation}
for any $X$, $Y$, $Z\in \Gamma(TM)$ and $V \in \Gamma(T^{\bot}M)$.

The analogue concept of locally product manifold is considered in the context of metallic geometry, having the name of \textit{locally metallic manifold} (\cite{Blaga2}). Thus, we say that the metallic Riemannian manifold $(\overline{M},\overline{g}, J)$ is \textit{locally metallic} if $J$ is parallel with respect to the Levi-Civita connection $\overline{\nabla}$ on $\overline{M}$ (i.e. $\overline{\nabla}J=0$).

\begin{remark}
In (\cite{Hr4}) we obtained that any almost product structure $F$ on $\overline{M}$ induces two metallic structures on $\overline{M}$:
\begin{equation}\label{e105}
J_{1}= \frac{2\sigma-p}{2}F+\frac{p}{2}I,
\end{equation}
\begin{equation}\label{e106}
J_{2}=-\frac{2\sigma-p}{2}F+\frac{p}{2}I,
\end{equation}
where $\sigma=\sigma _{p, q}=\frac{p+\sqrt{p^{2}+4q}}{2}$, with $p, q$ positive integer numbers.
\end{remark}

Also, for an almost product structure $F$ and for any $X \in \Gamma(TM)$ and $V \in \Gamma(T^{\perp}M)$,
the decompositions into the tangential and normal components of $FX$ and $FV$ are given by:
\begin{equation}\label{e107}
FX = fX + \omega X,
\end{equation}
\begin{equation}\label{e108}
FV = BV + CV,
\end{equation}
where $fX:=(F X)^T$, $\omega X:=(FX)^{\perp}$, $BV:=(F V)^T$ and $CV:=(F V)^{\perp}$.

Moreover, the maps $f$ and $C$ are $\overline{g}$-symmetric (\cite{Li&Liu}):
\begin{equation}\label{e109}
\overline{g}(fX,Y)=\overline{g}(X,fY),
\end{equation}
\begin{equation}
\overline{g}(CU,V)=\overline{g}(U,CV),
\end{equation}
for any $X, Y\in \Gamma(TM)$ and $U, V\in \Gamma(T^{\perp}M)$.

\begin{remark}(\cite{Hr6})
If $M$ is a submanifold in the almost product Riemannian manifold $(\overline{M}, \overline{g}, F)$ and $J$ is a metallic structure induced by $F$ on $\overline{M}$, then:
\begin{equation} \label{e110}
TX = \frac{p}{2}X \pm \frac{2\sigma-p}{2}fX,
\end{equation}
\begin{equation} \label{e111}
NX= \pm \frac{2\sigma-p}{2}\omega X,
\end{equation}
for any $X \in \Gamma(TM)$.
\end{remark}

\section{Pointwise slant submanifolds in metallic Riemannian manifolds}

B.-Y. Chen studied CR-submanifolds of a K\"{a}hler manifold which are warped products of holomorphic and totally real submanifolds, respectively (\cite{Chen3}, \cite{Chen1}, \cite{Chen2}). Also, in his new book (\cite{ChenBook}), he presents a multitude of properties for warped product manifolds and submanifolds, such as: warped product of Riemannian and K\"{a}hler manifolds, warped product submanifolds of K\"{a}hler manifolds (with the particular cases: warped product CR-submanifolds, warped product semi-slant or hemi-slant submanifolds of K\"{a}\-hler manifolds), CR-warped products in complex space forms and so on.

We shall state the notion of pointwise slant submanifold in a metallic Riemannian manifold,
following Chen's definition (\cite{Chen5}, \cite{Chen4}) of pointwise slant submanifold of an almost Hermitian manifold.

\begin{definition}
A submanifold $M$ of a metallic Riemannian manifold $(\overline{M}, \overline{g}, J)$ is called \textit{pointwise slant} if the angle $\theta_x(X)$ between $JX$ and $T_xM$ (called the \textit{Wirtinger angle}) is independent of the choice of the tangent vector $X \in T_{x}M\setminus\{0\}$, but it depends on $x \in M$. The Wirtinger angle is a real-valued function $\theta$ (called the Wirtinger function), verifying
\begin{equation}\label{e27}
\cos\theta_x =\frac{\| TX \|}{\| JX \|},
\end{equation}
for any $x\in M$ and $X \in T_{x}M\setminus\{0\}$.
\end{definition}

A pointwise slant submanifold of a metallic Riemannian manifold is called \textit{slant submanifold} if its Wirtinger function $\theta$ is globally constant.

\smallskip

In a similar manner as in (\cite{Chen5}) we obtain:
\begin{proposition}
If $M$ is an isometrically immersed submanifold in the metallic Riemannian manifold $(\overline{M}, \overline{g}, J)$, then $M$ is a pointwise slant submanifold if and only if
\begin{equation}\label{e28}
T^2=(\cos^2\theta)(pT+qI),
\end{equation}
for some real-valued function $\theta$.
\end{proposition}

From (\ref{e9}) and (\ref{e28}) we have:
\begin{proposition}
Let $M$ be an isometrically immersed submanifold in the metallic Riemannian manifold $(\overline{M}, \overline{g}, J)$. If $M$ is a pointwise slant submanifold with the Wirtinger angle $\theta$, then:
\begin{equation}\label{e29}
\overline{g}(NX,NY)=(\sin^2\theta)[p\overline{g}(TX,Y)+q\overline{g}(X,Y)]
\end{equation}
and
\begin{equation}\label{e30}
 tNX=(\sin^2\theta)(pTX+qX),
\end{equation}
for any $X$, $Y\in \Gamma(TM)$.
\end{proposition}

From (\ref{e28}), by a direct computation, we obtain:
\begin{proposition}
Let $M$ be an isometrically immersed submanifold in the metallic Riemannian manifold $(\overline{M}, \overline{g}, J)$. If $M$ is a pointwise slant submanifold with the Wirtinger angle $\theta$, then:
\begin{equation}\label{e31}
(\nabla_XT^2)Y=p(\cos^2\theta)(\nabla_XT)Y-\sin(2\theta)X(\theta)(pTY+qY),
\end{equation}
for any $X$, $Y\in \Gamma(TM)$.
\end{proposition}

\section{Pointwise bi-slant submanifolds in metallic Riemannian manifolds}

In this section we introduce the notion of pointwise bi-slant submanifold in the metallic context.

\begin{definition} \label{d2}
Let $M$ be an immersed submanifold in a metallic Riemannian manifold $(\overline{M},\overline{g},J)$. We say that $M$ is a {\it pointwise bi-slant submanifold} of $\overline{M}$ if there exists a pair of orthogonal distributions $D_{1}$ and $D_{2}$ on $M$ such that

(i) $TM = D_{1}\oplus D_{2}$;

(ii) $J(D_{1}) \bot D_{2}$ and $J(D_{2}) \bot D_{1}$;

(iii) the distributions $D_{1}$, $D_{2}$ are pointwise slant.
\end{definition}

If ${\theta_{1}}$ and ${\theta_{2}}$ are the slant functions of $D_{1}$ and $D_{2}$, respectively, then the pair
$\{ \theta_{1}, \theta_{2}\}$ is called the {\it bi-slant function}.

A pointwise slant submanifold $M$ is called {\it proper} if ${\theta_{1}}_x$, ${\theta_{2}}_x \neq 0; \frac{\pi}{2}$, for any $x\in M$ and both ${\theta_{1}}$, ${\theta_{2}}$ are not constant on $M$.

In particular, if ${\theta_{1}}=0$ and ${\theta_{2}} \neq 0; \frac{\pi}{2}$, then $M$ is called \textit{a pointwise semi-slant submanifold}; if ${\theta_{1}} = \frac{\pi}{2}$ and ${\theta_{2}}  \neq 0; \frac{\pi}{2}$, then $M$ is called \textit{a pointwise hemi-slant submanifold}.

\smallskip

Remark that if $M$ is a pointwise bi-slant submanifold of $\overline{M}$, then the distributions $D_{1}$ and $D_{2}$ on $M$ verify $T(D_{1}) \subseteq D_{1}$ and $T(D_{2}) \subseteq D_{2}$.

\begin{example}
Let $\mathbb{R}^{6}$ be the Euclidean space endowed with the usual Euclidean metric $\langle\cdot,\cdot\rangle$.
Let $i: M \rightarrow \mathbb{R}^{6}$ be the immersion given by:
$$i(u,v):=\left(\cos u \cos v, \cos u \sin v, \sin u \cos v, \sin u \sin v, \sin v, \cos v \right),$$
where $M :=\{(u, v) \mid  u, v \in (0, \frac{\pi}{2})\}$.

\pagebreak
A local orthogonal frame on $TM$ is given by:
 $$Z_{1}= -\sin u \cos v \frac{\partial}{\partial x_{1}} -\sin u \sin v  \frac{\partial}{\partial x_{2}}+ \cos u \cos v \frac{\partial}{\partial x_{3}} + \cos u \sin v  \frac{\partial}{\partial x_{4}}$$
 $$ Z_{2}= - \cos u \sin v \frac{\partial}{\partial x_{1}} +\cos u \cos v \frac{\partial}{\partial x_{2}} - \sin u \sin v \frac{\partial}{\partial x_{3}} +
 \sin u \cos v \frac{\partial}{\partial x_{4}} +$$$$+ \cos v \frac{\partial}{\partial x_{5}} - \sin v \frac{\partial}{\partial x_{6}}.$$

We define the metallic structure $J : \mathbb{R}^{6} \rightarrow \mathbb{R}^{6} $ by:
$$
 J(X_{1},X_{2},X_{3},X_{4},X_{5},X_{6}):=(\sigma X_{1}, \overline{\sigma} X_{2},\sigma X_{3},\overline{\sigma} X_{4},\sigma X_{5}, \overline{\sigma} X_{6} ),
 $$
where $\sigma:=\sigma_{p,q}=\frac{p+\sqrt{p^{2}+4q}}{2}$ is a metallic number ($p, q \in \mathbb{N}^{*}$) and $\overline{\sigma}=p-\sigma$.

We remark that $J$ verifies $J^{2}X=p J + q I$ and $\langle JX, Y\rangle = \langle X, JY\rangle$, for any $X$, $Y \in \mathbb{R}^{6}$.
Also, we have
 $$JZ_{1}=-\sigma\sin u \cos v \frac{\partial}{\partial x_{1}} -\overline{\sigma}\sin u \sin v  \frac{\partial}{\partial x_{2}}+ \sigma\cos u \cos v \frac{\partial}{\partial x_{3}} + \overline{\sigma}\cos u \sin v  \frac{\partial}{\partial x_{4}},$$
 $$ JZ_{2}= -\sigma \cos u \sin v \frac{\partial}{\partial x_{1}} +\overline{\sigma}\cos u \cos v \frac{\partial}{\partial x_{2}} - \sigma\sin u \sin v \frac{\partial}{\partial x_{3}} +$$ $$+\overline{\sigma}
 \sin u \cos v \frac{\partial}{\partial x_{4}} +\sigma \cos v \frac{\partial}{\partial x_{5}} -\overline{\sigma} \sin v \frac{\partial}{\partial x_{6}}.$$

 We remark that $\langle JZ_{1}, Z_{2}\rangle =\langle JZ_{2}, Z_{1}\rangle =0 $, $\langle JZ_{1}, Z_{1}\rangle= \sigma \cos^{2}v + \overline{\sigma}\sin^{2}v$ and $\langle JZ_{2}, Z_{2}\rangle = p$.

 On the other hand we get:
 $$\|Z_{1}\|=1, \ \ \|Z_{2}\|= \sqrt{2},$$
 $$\|J Z_{1}\|=\sqrt{\sigma^{2} \cos^{2}v + \overline{\sigma}^{2}\sin^{2}v}=\sqrt{p(\sigma \cos^{2}v + \overline{\sigma}\sin^{2}v)+q},$$ $$\|J Z_{2}\|=\sqrt{\sigma^{2}+\overline{\sigma}^{2}}=\sqrt{p^{2}+2q}.$$

We denote by $D_{1}:=span\{Z_{1}\}$ the pointwise slant distribution with the slant angle $\theta_{1}$, where $\cos \theta_{1} = \frac{f(u,v)}{\sqrt{pf(u,v)+q}}$, for $f(u,v):=\sigma \cos^{2}v + \overline{\sigma}\sin^{2}v$ a real function on $M$. Also, we denote by $D_{2}:=span\{Z_{2}\}$ the slant distribution with the slant angle $\theta_{2}$,
where $\cos \theta_{2} =\frac{p}{\sqrt{2(p^{2}+2q})}$.

The distributions $D_{1}$ and $D_{2}$ satisfy the conditions from Definition \ref{d2}.

If $M_{1}$ and $M_{2}$ are the integral manifolds of the distributions $D_{1}$ and $D_{2}$, respectively, then $M := M_{1} \times_{\sqrt{2}} M_{2}$ with the Riemannian metric tensor
$$g:= du^{2} +  2 d v^{2}$$ is a pointwise bi-slant submanifold in the metallic Riemannian manifold
$(\mathbb{R}^{6}, \langle\cdot,\cdot\rangle, J)$.
\end{example}

\begin{example}
In particular, if $f$ is a metallic function (i.e. $f^{2}=pf+q$), then $\cos \theta_{1}=1$ and we remark that $M$ is a semi-slant submanifold in the metallic Riemannian manifold $(\mathbb{R}^{6}, \langle\cdot,\cdot\rangle, J)$, with the slant angle $\theta=\theta_{2}$.
\end{example}

\begin{example}
On the other hand, if $f=0$ (i.e. $\tan v=\sqrt{- \frac{\sigma}{\overline{\sigma}}}=\frac{\sqrt{p^{2}+4q}+p}{2\sqrt{q}}$), then $\cos \theta_{1}=0$ and we remark that $M$ is a hemi-slant submanifold in the metallic Riemannian manifold $(\mathbb{R}^{6}, \langle\cdot,\cdot\rangle, J)$, with the slant angle $\theta=\theta_{2}$.
\end{example}

If we denote by $P_{i}$ the projections from $TM$ onto $D_{i}$, for $i \in \{1,2 \}$, then $X=P_{1}X+P_{2}X$, for any $X \in \Gamma(TM)$.
In particular, if $X \in D_{i}$, then $X=P_{i}X$, for $i \in \{1,2 \}$.

If we denote by $T_{i}=P_{i}\circ T$, for $i \in \{1,2 \}$, then, from (\ref{e4}), we obtain:
\begin{equation}\label{e32}
JX=T_{1}X+T_{2}X+NX.
\end{equation}

In a similar manner as in (\cite{Chen4}),
we get:
\begin{lemma}
Let $M$ be a pointwise bi-slant submanifold of a locally metallic Riemannian manifold $(\overline{M},\overline{g},J)$ with pointwise slant distributions $D_{1}$ and $D_{2}$, having slant functions ${\theta_{1}}$ and ${\theta_{2}}$. Then

(i) for any $X$, $Y\in D_1$ and $Z\in D_2$, we have:
\begin{equation}\label{e33}
(\sin^{2}\theta_{1}-\sin^{2}\theta_{2})\overline{g}(\nabla_{X}Y,pT_{2}Z+qZ)=
\end{equation}
$$=p[\overline{g}(\nabla_{X}Y,T_{2}Z)+\overline{g}(\nabla_{X}Z,T_{1}Y)]+
p(\cos^{2}\theta_{1}+1)\overline{g}(A_{NZ}Y+A_{NY}Z,X)-$$
$$-\overline{g}(A_{NT_{1}Y}Z+A_{NT_{2}Z}Y,X)-\overline{g}(A_{NZ}T_{1}Y+A_{NY}T_{2}Z,X);
$$

(ii) for any $X\in D_1$ and $Z$, $W\in D_2$, we have:
\begin{equation}\label{e34}
(\sin^{2}\theta_{2}-\sin^{2}\theta_{1})\overline{g}(\nabla_{Z}W,pT_{1}X+qX)=
\end{equation}
$$=p[\overline{g}(\nabla_{Z}W,T_{1}X)+\overline{g}(\nabla_{Z}X,T_{2}W)]+
p(\cos^{2}\theta_{2}+1)\overline{g}(A_{NX}W+A_{NW}X,Z)-$$
$$-\overline{g}(A_{NT_{2}W}X+A_{NT_{1}X}W,Z)-\overline{g}(A_{NW}T_{1}X+A_{NX}T_{2}W,Z).
$$
\end{lemma}

\proof
From (\ref{e1}) we have:
\begin{equation}\label{e35}
q\overline{g}(\nabla_{X}Y,Z)=q\overline{g}(\overline{\nabla}_{X}Y,Z)=\overline{g}(J^{2}\overline{\nabla}_{X}Y,Z)-p\overline{g}(J\overline{\nabla}_{X}Y,Z),
\end{equation}
for any $X$, $Y\in D_{1}$ and $Z \in D_{2}$.

By using (\ref{e2}) and $(\overline{\nabla}_{X}J)Y=0$, we obtain:
\begin{equation}\label{e36}
q\overline{g}(\nabla_{X}Y,Z)= \overline{g}(\overline{\nabla}_{X}J^2Y,Z)-p\overline{g}(\overline{\nabla}_{X}JY,Z).
\end{equation}

From (\ref{e32}) we get $JX = T_{1}X+NX$, $JY = T_{1}Y+NY$ and $JZ = T_{2}Z+NZ$, for any $X$, $Y \in D_1$ and $Z\in D_2$ and from here we obtain:
$$
q\overline{g}(\nabla_{X}Y,Z)=\overline{g}(\overline{\nabla}_{X}JT_{1}Y,Z)+\overline{g}(\overline{\nabla}_{X}JNY,Z) -p\overline{g}(\overline{\nabla}_{X}(T_{1}Y+NY),Z)=
$$
$$
=\overline{g}(\overline{\nabla}_{X}T_{1}^{2}Y,Z)+\overline{g}(\overline{\nabla}_{X}NT_{1}Y,Z)+\overline{g}(\overline{\nabla}_{X}NY,JZ)
-p\overline{g}(\overline{\nabla}_{X}T_{1}Y,Z)+$$$$+p\overline{g}(A_{NY}X,Z)=
$$
$$
=\overline{g}(\overline{\nabla}_{X}(\cos^{2} \theta_{1}(pT_{1}Y+qY)),Z)-\overline{g}(A_{NT_{1}Y}X,Z)+\overline{g}(\overline{\nabla}_{X}NY,T_{2}Z+NZ)+
$$
$$
+p\overline{g}(T_{1}Y,\overline{\nabla}_{X}Z)+p\overline{g}(A_{NY}X,Z).
$$

Thus, we get:
$$
q\overline{g}(\nabla_{X}Y,Z)=\cos^{2} \theta_{1}\overline{g}(\overline{\nabla}_{X}(pT_{1}Y+qY),Z)-\sin 2 \theta_{1}X(\theta_{1})\overline{g}(pT_{1}Y+qY,Z)-$$$$-\overline{g}(A_{NT_{1}Y}X,Z)-
\overline{g}(A_{NY}X,T_{2}Z)-\overline{g}(\overline{\nabla}_{X}NZ,JY)
+\overline{g}(\overline{\nabla}_{X}NZ,T_{1}Y)+$$$$+p\overline{g}(T_{1}Y,\overline{\nabla}_{X}Z)+p\overline{g}(A_{NY}X,Z).
$$

By using $\overline{g}(pT_{1}Y+qY,Z)=0$, we obtain:
$$
q\sin^{2}\theta_{1}\overline{g}(\nabla_{X}Y,Z)= p\cos^{2}\theta_{1}\overline{g}(\overline{\nabla}_{X}T_{1}Y,Z)-\overline{g}(A_{NT_{1}Y}Z+A_{NY}T_{2}Z,X)+
$$
$$
+\overline{g}(JNZ,\overline{\nabla}_{X}Y)-\overline{g}(A_{NZ}X,T_{1}Y)+p\overline{g}(T_{1}Y,\overline{\nabla}_{X}Z)+p\overline{g}(A_{NY}Z,X).
$$

Using (\ref{e10}) and (\ref{e30}), we find:
$$\overline{g}(JNZ,\overline{\nabla}_{X}Y) = \overline{g}(tNZ+nNZ,\overline{\nabla}_{X}Y) = $$
$$= \sin^{2}\theta_{2}\overline{g}(\nabla_{X}Y,qZ+pT_{2}Z)+\overline{g}(pNZ-NT_ 2Z,\overline{\nabla}_{X}Y)=$$
$$= q\sin^{2}\theta_{2}\overline{g}(\nabla_{X}Y,Z)+p\sin^{2}\theta_{2}\overline{g}(\nabla_{X}Y,T_{2}Z)-$$$$-p\overline{g}(\overline{\nabla}_{X}NZ,Y)+
\overline{g}(\overline{\nabla}_{X}NT_ 2Z,Y)$$
and from $$\overline{g}(\overline{\nabla}_{X}T_{1}Y,Z)=-\overline{g}(JY-NY,\overline{\nabla}_{X}Z)=\overline{g}(\overline{\nabla}_{X}Y,JZ)-\overline{g}(\overline{\nabla}_{X}NY,Z)=$$
$$=\overline{g}(\overline{\nabla}_{X}Y,T_ 2Z)-\overline{g}(Y,\overline{\nabla}_{X}NZ)-\overline{g}(\overline{\nabla}_{X}NY,Z) =$$
$$=\overline{g}(\overline{\nabla}_{X}Y,T_ 2Z)+\overline{g}(Y,A_{NZ}X)+\overline{g}(A_{NY}X,Z)
$$
we have:
$$
q(\sin^{2}\theta_{1}-\sin^{2}\theta_{2})\overline{g}(\nabla_{X}Y,Z)= p(1-\sin^{2}\theta_{1})\overline{g}(\nabla_{X}Y,T_{2}Z)+$$$$+p\cos^{2}\theta_{1}\overline{g}(A_{NZ}Y+A_{NY}Z,X)+\sin^{2}\theta_{2}\overline{g}(\nabla_{X}Y,pT_{2}Z)-
$$$$-\overline{g}(A_{NT_{1}Y}Z+A_{NY}T_{2}Z+A_{NZ}T_{1}Y+A_{NT_{2}Z}Y,X)-
$$
$$-p\overline{g}(Y,\overline{\nabla}_{X}NZ)+p\overline{g}(T_{1}Y,\nabla_{X}Z)+p\overline{g}(A_{NY}Z,X)$$
and from here we get (\ref{e33}).

In the same manner we find (\ref{e34}).

\begin{proposition}
Let $M$ be a pointwise semi-slant submanifold in a locally metallic Riemannian manifold $(\overline{M},\overline{g},J)$ with pointwise slant distributions $D_{1}$ and $D_{2}$, having slant functions $\theta_{1}$ and $\theta_{2}$.

(i) If $\theta_{1}=0$ and $\theta_{2}=\theta$, we obtain:
\begin{equation}\label{e350}
\sin^{2}\theta\overline{g}(\nabla_{X}Y,pT_{2}Z+qZ)=-p[\overline{g}(\nabla_{X}Y,T_{2}Z)+\overline{g}(\nabla_{X}Z,T_{1}Y)]
-
\end{equation}
$$
-2p\overline{g}(A_{NZ}Y,X)+\overline{g}(A_{NZ}T_{1}Y+A_{NT_{2}Z}Y,X),
$$
for any $X$, $Y\in D^T$ and $Z\in D^{\theta}$, and
\begin{equation}\label{e37}
\sin^{2}\theta\overline{g}(\nabla_{Z}W,pT_{1}X+qX)=p[\overline{g}(\nabla_{Z}W,T_{1}X)+\overline{g}(\nabla_{Z}X,T_{2}W)]+
\end{equation}
$$
+p(\cos^{2}\theta+1)\overline{g}(A_{NW}X,Z)-\overline{g}(A_{NT_{2}W}X+A_{NW}T_{1}X,Z),
$$
for any $X\in D^T$ and $Z$, $W \in D^{\theta}$.

(ii) If $\theta_{1}=\theta$ and $\theta_{2}=0$, we obtain:
\begin{equation}\label{e360}
\sin^{2}\theta\overline{g}(\nabla_{X}Y,T^{2}_{2}Z)=p\overline{g}(\nabla_{X}Y,T_{2}Z)-p\overline{g}(\nabla_{X}T_{1}Y,Z)+
\end{equation}
$$
+p(\cos^{2}\theta+1)\overline{g}(A_{NY}X,Z)-\overline{g}(A_{NT_{1}Y}Z+A_{NT_{2}Z}Y,X)-\overline{g}(A_{NY}T_{2}Z,X),
$$
for any $X, Y \in D^{\theta}$ and $Z\in D^{T}$, and
\begin{equation}\label{e370}
\sin^{2}\theta\overline{g}(\nabla_{Z}W,pT_{1}X+qX)=-p\overline{g}(\nabla_{Z}W,T_{1}X)+p\overline{g}(\nabla_{Z}T_{2}W,X)-
\end{equation}
$$
-2p\overline{g}(A_{NX}Z,W)+\overline{g}(A_{NT_{2}W}Z,X)+\overline{g}(A_{NT_{1}X}Z,W)+\overline{g}(A_{NX}T_{2}W,Z),
$$
for any $X\in D^{\theta}$ and $Z$, $W \in D^T$.
\end{proposition}

\begin{proposition}
Let $M$ be a pointwise hemi-slant submanifold in a locally metallic Riemannian manifold $(\overline{M},\overline{g},J)$ with pointwise slant distributions $D_{1}$ and $D_{2}$, having slant functions $\theta_{1}$ and $\theta_{2}$.

(i) If $\theta_{1}=\frac{\pi}{2}$ and $\theta_{2}=\theta$, we obtain:
\begin{equation}\label{e38}
\cos^{2}\theta\overline{g}(\nabla_{X}Y,pT_{2}Z+qZ)=p\overline{g}(\nabla_{X}Y,T_{2}Z)+
\end{equation}
$$+p\overline{g}(A_{NZ}Y+A_{NY}Z,X)-\overline{g}(A_{NT_{2}Z}Y+A_{NY}T_{2}Z,X),$$
for any $X$, $Y\in D^{\perp}$ and $Z\in D^{\theta}$, and
\begin{equation}\label{e39}
q\cos^{2}\theta\overline{g}(\nabla_{Z}W,X)=-p\overline{g}(\nabla_{Z}X,T_{2}W)-
\end{equation}
$$-
p(\cos^{2}\theta+1)\overline{g}(A_{NX}W+A_{NW}X,Z)+\overline{g}(A_{NT_{2}W}X+A_{NX}T_{2}W,Z),
$$
for any $X\in D^{\perp}$ and $Z$, $W\in D^{\theta}$.

(ii) If $\theta_{1}=\theta$ and $\theta_{2}=\frac{\pi}{2}$, we obtain:
\begin{equation}\label{e309}
q\cos^{2}\theta\overline{g}(\nabla_{X}Y,Z)=-p\overline{g}(\nabla_{X}Z,T_{1}Y)-
\end{equation}
$$-
p(\cos^{2}\theta+1)\overline{g}(A_{NZ}Y+A_{NY}Z,X)+\overline{g}(A_{NT_{1}Y}Z+A_{NZ}T_{1}Y,X),
$$
for any $X, Y\in D^{\theta}$ and $Z\in D^{\perp}$, and
\begin{equation}\label{e310}
\cos^{2}\theta\overline{g}(\nabla_{Z}W,pT_{1}X+qX)=p\overline{g}(\nabla_{Z}W,T_{1}X)+
\end{equation}
$$+p\overline{g}(A_{NX}W+A_{NW}X,Z)-\overline{g}(A_{NT_{1}X}W+A_{NW}T_{1}X,Z),$$
for any $X\in D^{\theta}$ and $Z$, $W\in D^{\perp}$.
\end{proposition}

\section{Warped product pointwise bi-slant submanifolds in metallic Riemannian manifolds}

In (\cite{Blaga1}), the authors of this paper introduced the Golden warped product Riemannian manifold and provided a necessary and sufficient condition for the warped product of two locally Golden Riemannian manifolds to be locally Golden. Moreover, the subject was continued in the papers (\cite{Blaga2}, \cite{Hr8}), where the authors characterized the metallic structure on the product of two metallic manifolds in terms of metallic maps and provided a necessary and sufficient condition for the warped product of two locally metallic Riemannian manifolds to be locally metallic.

Let $(M_1,g_1)$ and $({M_2},g_2)$ be two Riemannian manifolds (of dimensions $n_{1}>0$ and $n_{2}>0$, respectively) and let $\pi_1$,  $\pi_2$ be the projection maps from the product manifold ${M_1}\times {M_2}$ onto ${M_1}$ and ${M_2}$, respectively. We denote by $\widetilde{\varphi}:=\varphi \circ \pi_1$ the lift to ${M_1}\times {M_2}$ of a smooth function $\varphi$ on ${M_1}$. Then ${M_1}$ is called \textit{the base} and ${M_2}$ is called \textit{the fiber} of ${M_1}\times {M_2}$. The unique element $\widetilde{X}$ of $\Gamma(T({M_1}\times {M_2}))$ that is $\pi_1$-related to $X\in \Gamma(T{M_1})$ and to the zero vector field on ${M_2}$ will be called the \textit{horizontal lift of $X$} and the unique element $\widetilde{V}$ of $\Gamma(T({M_1}\times {M_2}))$ that is $\pi_2$-related to $V\in \Gamma(T{M_2})$ and to the zero vector field on ${M_1}$ will be called the \textit{vertical lift of $V$}.
We denote by $\mathcal{L}({M_1})$ the set of all horizontal lifts of vector fields on ${M_1}$ and by $\mathcal{L}({M_2})$ the set of all vertical lifts of vector fields on ${M_2}$.

For $f: M_1 \rightarrow (0,\infty)$ a smooth function on ${M_1}$, we consider the Riemannian metric $g$ on $M:={M_1}\times {M_2}$:
\begin{equation}\label{e40}
g:=\pi_1^* g_1+(f \circ \pi_1)^2 \pi_2^*g_2.
\end{equation}

\begin{definition} \label{d3}
The product manifold of ${M_1}$ and ${M_2}$ together with the Riemannian metric $g$ is called \textit{the warped product} of ${M_1}$ and ${M_2}$ by the warping function $f$ (\cite{Bishop}).

A warped product manifold $M:={M_1}\times_f {M_2}$ is called \textit{trivial} if the warping function $f$ is constant. In this case, $M$ is the Riemannian product ${M_1}\times {M_2}_f$, where ${M_2}_f$ is the manifold $M_2$ equipped with the metric $f^2 g_2$ (which is homothetic to $g_2$) (\cite{ChenBook}).
\end{definition}

In the next considerations, we shall denote by $(f \circ \pi_1)^2=: f^2$, $\pi_1^* g_1=: g_1$ and $\pi_2^*g_2=:g_2$, respectively.

\begin{lemma} (\cite{ChenBook})\label{1}
If $\nabla$ denotes the Levi-Civita connection on $M:={M_1}\times_f {M_2}$, then:
\begin{equation}\label{e41}
            \nabla_{X}Z=\nabla_{Z}X = X (\ln f) Z,
\end{equation}
for any $X$, $Y \in \Gamma(T{M_1})$ and $Z$, $W \in \Gamma(T{M_2})$.
\end{lemma}

The warped product $M_{1} \times_{f} M_{2}$ of two pointwise slant submanifolds $M_{1}$ and $M_{2}$ of a metallic Riemannian manifold $(\overline{M},\overline{g},J)$ is called a \textit{warped product pointwise bi-slant submanifold}. Moreover, it is called \textit{proper} if both $M_{1}$ and $M_{2}$ are proper pointwise slant submanifolds in $(\overline{M},\overline{g},J)$.

\smallskip

In a similar manner as in (\cite{Hr8}),
we get:
\begin{proposition}
Let $M:={M_1}\times_f {M_2}$ be a warped product pointwise bi-slant submanifold in a locally metallic Riemannian manifold $(\overline{M},\overline{g},J)$ with slant functions ${\theta_{1}}$, ${\theta_{2}}$ and warped function $f$. Then, for any $X$, $Y \in \Gamma(TM_{1})$ and $Z$, $W \in \Gamma(TM_{2})$, we have:
\begin{equation}\label{e42}
    \overline{g}(h(X,Y),NZ)=-\overline{g}(h(X,Z),NY),
\end{equation}
\begin{equation}\label{e43}
    \overline{g}(h(X,Z),NW)=0,
\end{equation}
\begin{equation}\label{e44}
     \overline{g}(h(Z,W),NX)= T_{1}X(\ln f)\overline{g}(Z,W)- X(\ln f)\overline{g}(Z,T_{2}W).
\end{equation}
\end{proposition}

\proof
For any $X$, $Y \in \Gamma(TM_{1})$ and $Z \in \Gamma(TM_{2})$, by using (\ref{e2}), (\ref{e4}), (\ref{e13}), (\ref{e41}) and $\overline{\nabla}J=0$ we obtain:
$$ \overline{g}(h(X,Y),NZ)= \overline{g}(\overline{\nabla}_{X}Y,JZ)-\overline{g}(\overline{\nabla}_{X}Y,T_{2}Z)=$$
$$=\overline{g}(\overline{\nabla}_{X}T_{1}Y,Z)+\overline{g}(\overline{\nabla}_{X}NY,Z)+\overline{g}(\overline{\nabla}_{X}T_{2}Z,Y)=$$
$$=-\overline{g}(\nabla_{X}Z,T_{1}Y)-\overline{g}(A_{NY}X,Z)+\overline{g}(Y,\nabla_{X}T_{2}Z)=$$
$$=-X(\ln f)\overline{g}(T_{1}Y,Z)-\overline{g}(h(X,Z),NY) + X(\ln f)\overline{g}(Y,T_{2}Z).$$

On the other hand, $\overline{g}(T_{1}Y,Z)=\overline{g}(JY,Z)=\overline{g}(Y,JZ)=\overline{g}(Y,T_{2}Z)$ and we obtain (\ref{e42}).

For any $X \in \Gamma(TM_{1})$ and $Z$, $W \in \Gamma(TM_{2})$, by using (\ref{e2}), (\ref{e4}), (\ref{e13}), (\ref{e41}) and $\overline{\nabla}J=0$ we obtain:
$$ \overline{g}(h(X,Z),NW)= \overline{g}(\overline{\nabla}_{X}Z,JW)-\overline{g}(\overline{\nabla}_{X}Z,T_{2}W)=$$
$$= \overline{g}(\nabla_{X}T_{2}Z,W)-\overline{g}(A_{NZ}X,W)-\overline{g}(\nabla_{X}Z,T_{2}W)=$$
$$= X (\ln f)[\overline{g}(T_{2}Z,W)-\overline{g}(Z,T_{2}W)]-\overline{g}(h(X,W),NZ)$$
and using
$$\overline{g}(T_{2}Z,W)-\overline{g}(Z,T_{2}W)=\overline{g}(JZ,W)-\overline{g}(Z,JW)=0,$$ we obtain
\begin{equation}\label{e500}
\overline{g}(h(X,Z),NW)= -\overline{g}(h(X,W),NZ).
\end{equation}

On the other hand, after interchanging $Z$ by $X$, we have:
 $$ \overline{g}(h(Z,X),NW)=\overline{g}(\nabla_{Z}T_{1}X,W)-\overline{g}(A_{NX}Z,W)-\overline{g}(\nabla_{Z}X,T_{2}W)=$$
 $$ =T_{1}X (\ln f)\overline{g}(Z,W)-X (\ln f)\overline{g}(Z,T_{2}W)-\overline{g}(h(Z,W),NX)=\overline{g}(h(X,W),NZ)$$
 and using (\ref{e500}) we get (\ref{e43}).

For any $X \in \Gamma(TM_{1})$ and $Z$, $W \in \Gamma(TM_{2})$, by using (\ref{e2}), (\ref{e4}), (\ref{e13}), (\ref{e41}) and $\overline{\nabla}J=0$ we obtain:
$$ \overline{g}(h(Z,W),NX)= \overline{g}(\overline{\nabla}_{Z}W,JX)-\overline{g}(\overline{\nabla}_{Z}W,T_{1}X)=$$
$$= \overline{g}(\nabla_{Z}T_{2}W,X)+\overline{g}(\overline{\nabla}_{Z}NW,X)-\overline{g}(\nabla_{Z}W,T_{1}X)=$$
$$= -\overline{g}(T_{2}W,\nabla_{Z}X)-\overline{g}(A_{NW}Z,X)+\overline{g}(W,\nabla_{Z} T_{1}X)=$$
$$ = - X(\ln f)\overline{g}(Z,T_{2}W)+T_{1}X(\ln f)\overline{g}(Z,W)$$
and we get (\ref{e44}).

\begin{proposition}
Let $M:={M_1}\times_f {M_2}$ be a warped product pointwise bi-slant submanifold in a locally metallic Riemannian manifold $(\overline{M},\overline{g},J)$ with slant functions ${\theta_{1}}$, ${\theta_{2}}$ and warped function $f$. Then, for any $X \in \Gamma(TM_{1})$ and $Z \in \Gamma(TM_{2})$, we have:
\begin{equation}\label{e45}
(\nabla_XT^2)Z=p(\cos^2\theta)(\nabla_XT)Z.
\end{equation}
\end{proposition}

\proof
From (\ref{e1}) we have:
\begin{equation}\label{e46}
q\overline{g}(\nabla_{X}Z,W)=q\overline{g}(\overline{\nabla}_{X}Z,W)=\overline{g}(J^{2}\overline{\nabla}_{X}Z,W)-p\overline{g}(J\overline{\nabla}_{X}Z,W),
\end{equation}
for any $X \in\Gamma(TM_{1})$ and $Z,W \in \Gamma(TM_{2})$.

By using (\ref{e2}), (\ref{e41}) and $(\overline{\nabla}_{X}J)Z=0$, we obtain:
\begin{equation}\label{e47}
q X(\ln f)\overline{g}(Z,W)= \overline{g}(\overline{\nabla}_{X}JZ,JW)-p\overline{g}(\overline{\nabla}_{X}JZ,W)
\end{equation}
and using $JZ=T_ 2 Z+ NZ$, for any $Z \in \Gamma(TM_{2})$, we have:
$$
q X(\ln f)\overline{g}(Z,W)= \overline{g}(\overline{\nabla}_{X}T_ 2 Z,T_ 2 W)+\overline{g}(\overline{\nabla}_{X}T_ 2Z,NW)+\overline{g}(\overline{\nabla}_{X}JNZ,W)-
$$
$$
-p\overline{g}(\overline{\nabla}_{X}T_ 2 Z,W)-p\overline{g}(\overline{\nabla}_{X}NZ,W).
$$

Thus, from (\ref{e10}) and (\ref{e41}) we get:
$$
q X(\ln f)\overline{g}(Z,W)=X(\ln f)\overline{g}(T_ 2Z,T_ 2W)+\overline{g}(h(X,T_ 2Z),NW)+$$$$+\overline{g}(\overline{\nabla}_{X}tNZ,W)+
\overline{g}(\overline{\nabla}_{X}nNZ,W)-pX(\ln f)\overline{g}(T_ 2Z,W)+p\overline{g}(A_{NZ}X,W).
$$

From (\ref{e43}) we obtain $\overline{g}(h(X,T_ 2Z), NW)=0$ and $\overline{g}(A_{NZ}X,W)=0$.

Thus, by using (\ref{e6}), (\ref{e28}) and (\ref{e30}), we have:
$$
q X(\ln f)\overline{g}(Z,W)=X(\ln f)\overline{g}(\cos^{2}\theta_ 2(pT_ 2Z+qZ),W )+$$
$$+\overline{g}(\overline{\nabla}_X(\sin^{2}\theta_ 2(pT_ 2Z+qZ)),W )+
 \overline{g}(\overline{\nabla}_X(pNZ-NT_ 2Z),W)-$$$$-pX(\ln f)\overline{g}(T_ 2Z,W)
$$
which implies
$$
\sin(2\theta_ 2)X(\theta_ 2)\overline{g}(pT_ 2Z+qZ,W )=p\overline{g}(h(X,W),NZ)-\overline{g}(h(X,W),T_2Z)=0.
$$
Thus, from (\ref{e43}) and (\ref{e31}) we get (\ref{e45}).

\begin{example}
Let $\mathbb{R}^{6}$ be the Euclidean space endowed with the usual Euclidean metric $\langle\cdot,\cdot\rangle$.
Let $i: M \rightarrow \mathbb{R}^{6}$ be the immersion given by:
$$i(u,v):=\left( u \sin v, u \cos v, u, u \cos v, u \sin v, v \right),$$
where $M :=\{(u, v) \mid  u>0, v \in (0, \frac{\pi}{2})\}$.

A local orthogonal frame on $TM$ is given by:
 $$Z_{1}= \sin v \frac{\partial}{\partial x_{1}} + \cos v  \frac{\partial}{\partial x_{2}}+ \frac{\partial}{\partial x_{3}} + \cos v  \frac{\partial}{\partial x_{4}}+\sin v  \frac{\partial}{\partial x_{5}}$$
 $$ Z_{2}= u \cos v \frac{\partial}{\partial x_{1}} - u \sin v \frac{\partial}{\partial x_{2}} - u \sin v \frac{\partial}{\partial x_{4}} + u \cos v \frac{\partial}{\partial x_{5}}+\frac{\partial}{\partial x_{6}}.$$

We define the metallic structure $J : \mathbb{R}^{6} \rightarrow \mathbb{R}^{6} $ by:
$$
 J(X_{1},X_{2},X_{3},X_{4},X_{5},X_{6}):=(\sigma X_{1}, \sigma X_{2},\sigma X_{3},\overline{\sigma} X_{4}, \overline{\sigma}X_{5}, \overline{\sigma} X_{6} ), $$
where $\sigma:=\sigma_{p,q}=\frac{p+\sqrt{p^{2}+4q}}{2}$ is a metallic number ($p, q \in \mathbb{N}^{*}$) and $\overline{\sigma}=p-\sigma$.

We remark that $J$ verifies $J^{2}X=p J + q I$ and $\langle JX, Y\rangle = \langle X, JY\rangle$, for any $X$, $Y \in \mathbb{R}^{6}$.
Also, we have:
 $$J Z_{1}= \sigma \sin v \frac{\partial}{\partial x_{1}} +\sigma \cos v  \frac{\partial}{\partial x_{2}}+ \sigma \frac{\partial}{\partial x_{3}} + \overline{\sigma} \cos v  \frac{\partial}{\partial x_{4}}+
 \overline{\sigma} \sin v  \frac{\partial}{\partial x_{5}}$$
 $$ J Z_{2}=  \sigma u\cos v \frac{\partial}{\partial x_{1}} - \sigma u\sin v \frac{\partial}{\partial x_{2}} - \overline{\sigma} u\sin v \frac{\partial}{\partial x_{4}} +  \overline{\sigma} u \cos v \frac{\partial}{\partial x_{5}} + \overline{\sigma} \frac{\partial}{\partial x_{6}}.$$

 We remark that $\langle JZ_{1}, Z_{2}\rangle =\langle JZ_{2}, Z_{1}\rangle =0 $, $\langle JZ_{1}, Z_{1}\rangle= 2\sigma +\overline{\sigma}$ and
 $\langle JZ_{2}, Z_{2}\rangle = u^{2}(\sigma+\overline{\sigma})+\overline{\sigma}$.

 On the other hand we get:
  $$\|Z_{1}\|=\sqrt{3}, \ \ \|Z_{2}\|= \sqrt{2u^{2}+1},$$
   $$\|J Z_{1}\|=\sqrt{2 \sigma^{2} + \overline{\sigma}^{2}}, \ \ \|J Z_{2}\|=\sqrt{u^{2}(\sigma^{2}+\overline{\sigma}^{2})+\overline{\sigma}^{2}}.$$

We denote by $D_{1}:=span\{Z_{1}\}$ the slant distribution with the slant angle $\theta_{1}$, where $\cos \theta_{1} = \frac{2\sigma+\overline{\sigma}}{\sqrt{3(2 \sigma^{2} + \overline{\sigma}^{2})}}$. Also, we denote by $D_{2}:=span\{Z_{2}\}$ the pointwise slant distribution with the slant angle $\theta_{2}$, where
$\cos \theta_{2} = \frac{u^{2}(\sigma+\overline{\sigma})+\overline{\sigma}}{\sqrt{(2u^{2}+1)(u^{2}(\sigma^{2}+\overline{\sigma}^{2})+\overline{\sigma}^{2})}}$.

The distributions $D_{1}$ and $D_{2}$ satisfy the conditions from Definition \ref{d3}.

If $M_{1}$ and $M_{2}$ are the integral manifolds of the distributions $D_{1}$ and $D_{2}$, respectively, then $M := M_{1} \times_{\sqrt{2u^{2}+1}} M_{2}$ with the Riemannian metric tensor
$$g:= 3 du^{2} + (2u^{2}+1) d v^{2}$$ is a warped product pointwise bi-slant submanifold in the metallic Riemannian manifold $(\mathbb{R}^{6}, \langle\cdot,\cdot\rangle, J)$.
\end{example}

\section{Warped product pointwise semi-slant or hemi-slant submanifolds in metallic Riemannian manifolds}

In this section we get some properties of
pointwise semi-slant and pointwise hemi-slant submanifolds in locally metallic Riemannian manifolds.

\begin{definition} \label{d30}
Let $M:={M_1}\times_f {M_2}$ be a warped product bi-slant submanifold in a metallic Riemannian manifold $(\overline{M},\overline{g},J)$ such that one of the components $M_{i}$ ($i \in \{1,2\}$) is an invariant submanifold (respectively, anti-invariant submanifold) in $\overline{M}$ and the other one is a pointwise slant submanifold in $\overline{M}$, with the Wirtinger angle $\theta_x \in [0, \frac{\pi}{2}]$. Then we call the submanifold $M$ \textit{warped product pointwise semi-slant submanifold} (respectively, \textit{warped product pointwise hemi-slant submanifold}) in the metallic Riemannian manifold $(\overline{M},\overline{g},J)$.
\end{definition}

In a similar manner as in Theorem 2 from (\cite{Hr8}), we obtain:

\begin{theorem}
If $M:={M_{T}}\times_f {M_{\theta}}$ is a warped product pointwise semi-slant submanifold in a locally metallic Riemannian manifold $(\overline{M},\overline{g},J)$ with the pointwise slant angle $\theta_ x \in (0, \frac{\pi}{2})$, for $x \in M_{\theta}$, then the warping function $f$ is constant on the connected components of $M_{T}$.
\end{theorem}

\proof
For any $X \in \Gamma(TM_{T})$, $Z \in \Gamma(TM_{\theta})\setminus\{0\}$, by using (\ref{e13}) in $\overline{\nabla}_{Z}JX = J\overline{\nabla}_{Z}X$ and (\ref{e41}), we obtain:
$$
TX(\ln f)Z+h(TX,Z)=T\nabla_{Z}X+N\nabla_{Z}X+t h(X,Z)+n h(X,Z).
$$

From the equality of the normal components of the last equation, it follows
\begin{equation}\label{e320}
 h(TX,Z) = X(\ln f)NZ + n h(X,Z)
\end{equation}
and replacing $X$ with $TX=JX$ (for $X \in \Gamma(TM_{T})$) in (\ref{e320}), we obtain:
 $$
 h(J^{2}X,Z)=TX(\ln f)NZ+n h(TX,Z).
 $$

 Thus, we get:
 $$
TX(\ln f)\overline{g}(NZ,NZ)= \overline{g}(h(J^{2}X,Z),NZ)-\overline{g}(n h(TX,Z),NZ)=
 $$
 $$
=p \overline{g}(h(TX,Z),NZ) +q \overline{g} (h(X,Z),NZ) -\overline{g}(nh(TX,Z),NZ),
 $$
 for any $X \in \Gamma(TM_{T})$ and $Z \in \Gamma(TM_{\theta})$.

From (\ref{e43}) we have $\overline{g}(h(TX,Z),NZ)= \overline{g} (h(X,Z),NZ) =0$, for any $X \in \Gamma(TM_{T})$ and $Z \in \Gamma(TM_{\theta})$ and by using (\ref{e29}), we get:
 \begin{equation}\label{e303}
TX(\ln f)\sin^{2}\theta [p \overline{g}(TZ,Z)+q \overline{g}(Z,Z)] = -\overline{g}(nh(TX,Z),NZ).
 \end{equation}

On the other hand, for any $X \in \Gamma(TM_{T})$ and $Z \in \Gamma(TM_{\theta})$, we have $TX \in \Gamma(TM_{T})$ and $TZ \in \Gamma(TM_{\theta})$ and from (\ref{e43}), we obtain:
$$\overline{g}(h(TX,Z),NZ)=\overline{g}(h(TX,Z),NTZ)=0.$$

Thus, by using (\ref{e1}) and (\ref{e7}), we have:
 $$\overline{g}(nh(TX,Z),NZ)=\overline{g}(h(TX,Z),nNZ)=\overline{g}(h(TX,Z),J^{2}Z-JTZ)=$$
 $$=p\overline{g}(h(TX,Z),NZ)+q\overline{g}(h(TX,Z),Z)-\overline{g}(h(TX,Z),NTZ)=0$$
 and using (\ref{e303}), we obtain:
 $$TX(\ln f)\tan^{2}\theta_x  \overline{g}(TZ,TZ)=0,$$ for any $Z \in \Gamma(TM_{\theta})$ and $x \in M_{\theta}$.

Since $\theta_x \in (0, \frac{\pi}{2})$ and $TZ \neq  0$, we get $TX(\ln f) = 0$, for any $X \in \Gamma(TM_{T})$, which implies that the warping function $f$ is constant on the connected components of $M_{T}$.

\begin{theorem}
If $M:={M_\theta}\times_f {M_T}$ is a warped product pointwise semi-slant submanifold in a locally metallic Riemannian manifold $(\overline{M},\overline{g},J)$ with the pointwise slant angle $\theta_x \in (0, \frac{\pi}{2})$, for $x \in M_\theta$, then
$$(A_{NT_{1}Y}X - A_{NT_{1}X}Y) \in \Gamma(TM_{\theta}),$$ for any $X,Y \in \Gamma(TM_{\theta})$.
\end{theorem}

\proof
For any $X, Y\in \Gamma(TM_{\theta})$ and $Z\in \Gamma(TM_{T})\setminus\{0\}$, from (\ref{e36}) and the symmetry of the shape operator, we have:
$$
\sin^{2}\theta\overline{g}([X,Y],T^{2}_{2}Z)=p\overline{g}([X,Y],T_{2}Z)-p\overline{g}(\nabla_{X}T_{1}Y-\nabla_{Y}T_{1}X,Z)
+
$$
$$
+p(\cos^{2}\theta+1)[\overline{g}(h(X,Z),NY)-\overline{g}(h(Y,Z),NX)]-\overline{g}(h(X,Z),NT_{1}Y)+
$$
$$+\overline{g}(h(Y,Z),NT_{1}X)+\overline{g}(h(X,Y),NT_{2}Z)-\overline{g}(h(Y,X),NT_{2}Z)-$$
$$-\overline{g}(h(X,T_{2}Z),NY)+\overline{g}(h(Y,T_{2}Z),NX).$$

Using (\ref{e2}) and (\ref{e42}), we obtain:
$$
\overline{g}(\nabla_{X}T_{1}Y-\nabla_{Y}T_{1}X,Z)=\overline{g}(\nabla_{X}JY-\nabla_{Y}NX-\nabla_{Y}JX+\nabla_{Y}NX,Z)=
$$
$$
=\overline{g}(\nabla_{X}Y,JZ)-\overline{g}(\nabla_{Y}X,JZ)+\overline{g}(A_{NY}X,Z)-\overline{g}(A_{NX}Y,Z)=
$$
$$
=\overline{g}([X,Y],JZ)+\overline{g}(h(X,Z),NY)-\overline{g}(h(Z,Y),NX)=\overline{g}([X,Y],T_{2}Z).
$$

From (\ref{e42}) we get: $$\overline{g}(h(X,Z),NY)=\overline{g}(h(Y,Z),NX)=-\overline{g}(h(X,Y),NZ).$$

Thus, using the symmetry of the shape operator, we have:
$$
\overline{g}(h(X,T_{2}Z),NY)-\overline{g}(h(Y,T_{2}Z),NX)=$$$$=-\overline{g}(h(X,Y)NT_{2}Z)+\overline{g}(h(Y,X),NT_{2}Z)=0
$$
and
$$\overline{g}(h(X,Z),NT_{1}Y)-\overline{g}(h(Y,Z),NT_{1}X)=\overline{g}(A_{NT_{1}Y}X-A_{NT_{1}X}Y,Z).$$

Thus, we obtain:
 $$
\sin^{2}\theta\overline{g}([X,Y],T^{2}_{2}Z)=\overline{g}(A_{NT_{1}Y}X-A_{NT_{1}X}Y,Z),
$$
which implies the conclusion.

Following the same steps such in (\cite{Hr8}), we can prove that:

\begin{theorem}
If $M:={M_{\perp}}\times_f {M_{\theta}}$ (or $M:={M_{\theta}}\times_f {M_{\perp}}$) is a warped product pointwise hemi-slant submanifold in a locally metallic Riemannian manifold $(\overline{M},\overline{g},J)$ with the pointwise slant angle $\theta_x \in (0, \frac{\pi}{2})$, for $x \in M_\theta$, then the warping function $f$ is constant on the connected components of $M_{\perp}$ if and only if
\begin{equation}\label{e10008}
A_{NZ}X = A_{NX}Z,
\end{equation}
for any $X \in \Gamma(TM_{\perp})$ and $Z \in \Gamma(TM_{\theta})$ (or $X \in \Gamma(TM_{\theta})$ and $Z \in \Gamma(TM_{\perp})$, respectively).
\end{theorem}

\bigskip

\textit{Cristina E. Hretcanu}

\textit{Stefan cel Mare University of Suceava, Romania}

\textit{criselenab@yahoo.com}

\bigskip

\textit{Adara M. Blaga}


\textit{West University of Timi\c{s}oara, Rom\^{a}nia}

\textit{adarablaga@yahoo.com}

\end{document}